\documentclass[12pt]{article}
\usepackage{amsmath}
\usepackage{amssymb}
\usepackage{graphicx}
\usepackage[cp1251]{inputenc}
\usepackage[russian]{babel}
\newcommand{\il}[2]{\int\limits_{#1}^{#2}}
\newcommand{\ilp}[1]{\int\limits_{#1}^{+\infty}}

\newcommand{\ph}{\phantom{a}}
\newcommand{\phh}{\phantom{aaa}}

\newcommand{\sist}[2]{\left\{
\begin{array}{l}
{#1}\\
\ph\\
{#2}
\end{array}
\right.}

\begin{document}

\vskip 20pt

MSC 34C10

\vskip 10pt

\centerline{\bf On the oscillation  of}
 \centerline{\bf linear matrix Hamiltonian systems}

\vskip 10 pt

\centerline{\bf G. A. Grigorian}
\centerline{\it Institute  of Mathematics NAS of Armenia}
\centerline{\it E -mail: mathphys2@instmath.sci.am}
\vskip 10 pt

\noindent
Abstract. The Riccati equation method is used to establish new oscillation criteria for linear matrix Hamiltonian systems. New approaches allow to extend and completed a result, obtained by S. Kumary and S. Umamaheswaram. The oscillation problem for linear matrix Hamiltonian systems in a new direction, which is to break the positive definiteness condition, imposed on one of the coefficients of the system, is investigated.
Some examples are provided for comparing the obtained results with each other and with the result of S. Kumary and S. Umamaheswaram, as well as to illustrate the applicability of these results.

\vskip 10 pt


\vskip 10 pt

{\bf 1. Introduction}.
Let $A(t), \ph B(t)$ and $C(t)$ be complex-valued locally integrable   matrix  functions of dimension $n \times n$ on $[t_0,+\infty)$ and let $B(t)$ and $C(t)$ be Hermitian,  i.e.  $B(t) = B^*(t), \ph C(t) = C^*(t), \ph t\ge t_0$, where $*$ denotes the conjugation sign. Consider the linear matrix Hamiltonian system
$$
\sist{\Phi' = A(t) \Phi + B(t) \Psi;}{\Psi' = C(t)\Phi - A^*(t) \Psi, \phh t\ge t_0.} \eqno (1.1)
$$
 By a solution $(\Phi(t), \Psi(t))$ of this system we mean a pair of absolutely continuous matrix functions $\Phi = \Phi(t)$ and $\Psi = \Psi(t)$ of dimension $n\times n$ on $[t_0, +\infty)$ satisfying (1.1) almost everywhere on $[t_0,+ \infty)$.

{\bf Definition 1.1}. {\it A solution $(\Phi(t), \Psi(t))$ of the system (1.1) is called conjoined (or prepared, preferred), if $\Phi^*(t) \Psi(t) = \Psi^*(t) \Phi(t), \ph t\ge t_0.$}

{\bf Definition 1.2}. {\it The system (1.1) is called oscillatory if  for its every conjoined solution $(\Phi(t), \Psi(t))$  the function $\det \Phi(t)$ has arbitrary large zeroes.}

The oscillation problem for linear matrix Hamiltonian systems is that of finding explicit conditions on the coefficients of the system, providing its oscillation. This is an important problem
of qualitative theory of differential equations and many works are devoted to it (see e.g., [1, 3, 5, 6, 13-18] and cited works therein).  Among them notice the following result, obtained by S. Kumari   and S.Umamaheswaram

{\bf Theorem 1.1 ([11, Theorem 2.1])}. {\it Let $A(t), \ph B(t)$ and $C(t)$ be real-valued continuous functions on $[t_0,+\infty)$ and let $B(t)$ be positive definite for all $t \ge t_0$. If there exists a positive linear functional $g$ on the space of matrices of dimension $n\times n$ such that
$$
\lim\limits_{t \to +\infty}\il{t_0}{t}\frac{d s}{g[B^{-1}(s)]} = +\infty,
$$
$$
\lim\limits_{t\to +\infty}g\left[ -\il{t_0}{t} \Bigl(C(s) + A^*(s) B^{-1}(s) A(s)\Bigr) ds - B^{-1}(t) A(t)\right] = +\infty,
$$
then the system (1.1) is oscillatory.
}

In this paper we use the Riccati equation method for obtaining some new oscillation criteria for the system (1.1). We show that one of the obtained results is an extension of Theorem 1.1.
Traditionally the oscillation problem for the system (1.1) was studied under the restriction that the coefficient $B(t)$ of the system (1.1) is positive definite (a condition, which is essential from the point of view if used methods). New approaches in the papers [5] and [6]  allowed to obtain oscillation criteria in the direction of weakening (braking) the positive definiteness of $B(t)$. In this paper we continue to study the oscillation problem for the system (1.1) in the mentioned direction. On examples we compare the obtained oscillation criteria with each other and with Theorem 1.1 and demonstrate their applicability.

{\bf 2. Auxiliary propositions.} Let $a(t), \ph b(t), \ph c(t), \ph  a_1(t), \ph b_1(t)$ and $c_1(t)$ be real-valued locally integrable functions on $[t_0,+\infty)$. Consider the Riccati equations.
$$
y' + a(t) y^2 + b(t) y + c(t) = 0, \phh t \ge t_0,  \eqno (2.1)
$$
$$
y' + a_1(t)y^2 + b_1(t) y + c_1(t) = 0, \phh t \ge t_0 \eqno (2.2)
$$
and the differential inequalities
$$
\eta' + a(t) \eta^2 + b(t) \eta + c(t) \ge 0, \phh t \ge t_0,  \eqno (2.3)
$$
$$
\eta' + a_1(t)\eta^2 + b_1(t) \eta + c_1(t) \ge 0, \phh t \ge t_0 \eqno (2.4)
$$
Let $[t_1,t_2)$ be an interval in $[t_0,+\infty) \ph (t_0 \le t_1 < t_2 \le +\infty)$. By a solution of Eq. (2.1) (of Eq. (2.2), inequality (2.3), inequality (2.4)) on $[t_1,t_2)$  we mean an absolutely continuous function on $[t_1,t_2)$, satisfying (2.1) ((2.2), (2,3), (2.4)) almost everywhere on $[t_1,t_2)$. Note that every solution of Eq. (2.1) (Eq. (2.2)) on any interval $[t_1,t_2)$ is also a solution of the inequality (2.3) ((2.4)) on that interval. Note also that for ($a(t) \ge 0 \ph (a_1(t) \ge 0), \ph t\ge t_0$ the real-valued solutions of the linear equation
$$
\eta' + b(t) \eta + c(t) = 0, \ph t \ge t_0 \phh  (\eta' + b_1(t) \eta + c_1(t) = 0, \ph t \ge t_0)
$$
are solutions of the inequality (2.3) ((2.4)) as well. Therefore, for $a(t)\ge 0, \ph (a_1(t) \ge 0), \ph t \ge t_0$
the inequality (2.3) ((2.4)) has a solution on $[t_0,+\infty)$, satisfying any initial-value condition.

{\bf Theorem 2.1.} {\it Let $y_0(t)$ be a solution of Eq. (2.1) on $[t_1,t_2)$, $\eta_0(t), \ph \eta_1(t)$ be a solution of the inequalities (2.3) and (2.4) respectively with $\eta_0(t_1) \ge y_0(t_1), \ph \eta_1(t_1) \ge y_0(t_1)$, and let $a_1(t) \ge0,\ph \lambda - y_0(t_1) + \il{t_1}{t}\exp\biggl\{\il{t_1}{\tau}[a_1(\xi)(\eta_0(\xi) + \eta_1(\xi)) = b_1(\xi)]d\xi\biggr\}\times$\\
$\times [(a(\tau) - a_1(\tau)) y^2_0(\tau) + (b(\tau)- b_1(\tau)) y_0(\tau) = c(\tau) - c_1(\tau)] d \tau \ge 0, \ph t \in [t_1,t_2)$, for some $\lambda \in [y_0(t_1),\eta_1(t_1)].$ Then Eq. (2.2) has a solution $y_1(t)$ on $[t_1,t_2)$ with $y_1(t_1) \ge y_0(t_1)$, moreover $y_1(t) \ge y_0(t), \ph t \in [t_1,t_2)$.}

Proof. By analogy with the proof of Theorem 3.1 from [8].

Let $a_{jk}(t), \ph j,k =1,2$ be real-valued locally integrable functions on $[t_0,+\infty)$. Consider the linear system of ordinary differential equations
$$
\sist{\phi' = a_{11}(t) \phi + a_{12}(t) \psi,}{\psi' = a_{21}(t) \phi + a_{22}(t) \psi, \ph t\ge t_0.} \eqno (2.5)
$$
and the corresponding Riccati equation
$$
y' + a_{12}(t) y^2 + E(t) y -a_{21}(t) = 0, \phh t \ge t_0, \eqno (2.6)
$$
where $E(t) \equiv a_{11}(t) - a_{22}(t), \ph t \ge t_0$. By a solution of the system (2.5) we mean an ordered pair $(\phi(t), \psi(t))$ of absolutely continuous functions $\phi(t), \ph \psi(t)$ on $[t_0,+\infty)$, satisfying for $\phi =\phi(t), \ph \psi = \psi(t)$  the system (2.5) almost everywhere on $[t_0,+\infty)$.
All solutions $y(t)$ of Eq. (2.6), existing on any interval $[t_1,t_2) \subset [t_0,+\infty)$ are connected with solutions $(\phi(t),\psi(t))$  of the system (2.5) by relations (see [8])
$$
\phi(t) = \phi(t_1)\exp\biggl\{\il{t_1}{t}[a_{12}(\tau) y(\tau) + a_{11}(\tau)]d\tau\biggr\}, \ph \phi(t_1) \ne 0, \phh  \psi(t) = y(t) \phi(t),  \eqno (2.7)
$$
$t \in [t_1,t_2).$

{\bf Definition 2.1.} {\it The system (2.2) is called oscillatory if for its  every solution  $(\phi(t), \psi(t))$ the function $\phi(t)$ has arbitrary large zeroes.}

\vskip 5pt

{\bf Remark 2.3}. {\it Explicit oscillatory criteria for the system (2.1) (therefore for the system (2.5)) are obtained in [7].}

{\bf Theorem 2.2.} {\it Let the following conditions be satisfied.

\noindent
$a(t) \ge 0, \ph t \ge t_0.$
$$
\ilp{t_0}a_{12}(t)a_{12}(t)\exp\biggl\{-\il{a}{t}E(\tau) d \tau\biggr\} d t = - \ilp{t_0}a_{21}(t)\exp\biggl\{\il{a}{t}E(\tau) d \tau\biggr\}d t = +\infty.
$$
}

Proof. By analogy with the proof of Corollary 3.1 from [10] (see also [6, Theorem 2.4]).

{\bf Definition 2.2.} {\it An interval $[t_1,t_2)\subset [t_0,+\infty)$ is called the maximum existence interval for a solution $y(t)$ of Eq. (2.1), if $y(t)$ exists on $[t_1,t_2)$ and cannot be continued to the right from $t_2$ as a solution of Eq. (2.1).
}

{\bf Lemma 2.1.} {\it Let $y(t)$ be a solution of Eq. (2.1) on $[t_1,t_2) \subset [t_0,+\infty)$, and let $t_2< +\infty$. If  $a(t) \ge 0, \ph t \in [t_1,t_2)$ and  the function $F(t) \equiv \il{t_1}{t}a(\tau) y(\tau) d\tau, \ph t \in [t_1,t_2)$ is bounded from below on $[t_1,t_2)$, then $[t_1,t_2)$ is not the maximum existence interval for $y(t)$.
}

Proof. By analogy with the proof of Lemma 2.1 from [9].

Let $e(t)$ and $e_1(t)$ be real-valued functions on $[t_0,+\infty)$ and let $e(t)$ be locally integrable and $e_1(t)$ be absolutely continuous on $[t_0,+\infty)$. Consider the Riccati integral equations
$$
y(t) + \il{t_0}{t} a(\tau) y^2(\tau) d\tau + e(t) = 0, \phh t \ge t_0, \eqno (2.8)
$$
$$
y(t) + \il{t_0}{t} a(\tau) y^2(\tau) d\tau + e_1(t) = 0, \phh t \ge t_0, \eqno (2.9)
$$

{\bf Lemma 2.2.} {\it Let $y_0(t)$ be a solution of Eq. (2.8) on $[t_0,t_1)$. If $a(t) \ge 0, \ph e(t) > e_1(t) > 0, \ph t\in[t_0,t_1)$, then Eq. (2.9) has a solution $y_1(t)$  on $[t_0,t_1)$ and
$$
y_1(t) > y_0(t), \phh t\in [t_0,t_1). \eqno (2.10)
$$
}

Proof. Since $a(t) \ge 0, \ph e(t) > 0, \ph t\in [t_0,t_1)$ by (2.8)
$$
y_0(t) < 0, \phh t \in [t_0,t_1). \eqno (2.11)
$$
Let $y_1(t)$ be a solution of Eq. (2.9) and let $[t_0,t_2)$ be its maximum existence interval. Show that
$$
t_2 \ge t_1. \eqno (2.12)
$$
Suppose
$$
t_2 < t_1. \eqno (2.13)
$$
Show that
$$
y_1(t) > y_0(t), \ph t\in [t_0,t_2). \eqno (2.14)
$$
Suppose that this is false. By (2.8) and (2.9) from the conditions $e(t) > e_1(t) > 0, \linebreak t\in [t_0,t_1)$ of the lemma it follows that $y_1(t_0) = - e_1(t_0) > - e(t_0) = y_0(t_0)$. Then there exists $t_3\in(t_0,t_2)$ such that
$$
y_1(t) > y_0(t), \ph t_0 \le t < t_3, \eqno (2.15)
$$
$$
y_1(t_3) = y_0(t_3). \eqno (2.16)
$$
On the other hand by (2.8) and (2.9) we have $y_1(t_3) - y_0(t_3) = \il{t_0}{t_3}a(\tau)[y_0^2(\tau) - y_1^2(\tau)]d\tau + e(t_3) - e_1(t_3)$. This together with (2.11),     (2.15) and the condition $a(t) \ge 0, \ph t \in [t_0,t_1)$  of the lemma implies that $y_1(t_3) > y_0(t_1)$, which contradicts (2.16). The obtained contra-\linebreak diction proves (2.14). Obviously $y_1(t)$ is a solution of the Riccati equation (recall that $e_1(t)$ is absolutely continuous)
$$
y' + a(t) y^2 + e'_1(t) = 0, \phh t \ge t_0
$$
on $[t_0,t_1)$. Then by Lemma 2.1  and (2.7) from the condition $a(t)\ge 0, \ph t\in [t_0,t_1)$ of the lemma and from (2.6) it follows that $[t_0,t_2)$ is not the maximum existence interval for $y_1(t)$, which contradicts our supposition. The obtained contradiction proves (2.12). From (2.12) and (2.14) it follows existence $y_1(t)$ on $[t_0,t_1)$ and the inequality (2.10). The lemma is proved.

{\bf Lemma 2.3.} {\it  For any two square matrices  $M_k \equiv (m_{ij}^l)_{i,j=1}^n, \ph l=1,2$  he equality
$$
tr(M_1 M_1) = tr (M_2 M_1)
$$
is valid.}

Proof.  We have $tr (M_1 M_2) = \sum\limits_{j=1}^n(\sum\limits_{k=1}^n m_{jk}^1 m_{kj}^2) = \sum\limits_{k=1}^n(\sum\limits_{j=1}^n m_{jk}^1 m_{kj}^2) = \sum\limits_{k=1}^n(\sum\limits_{j=1}^n m_{kj}^2 m_{jk}^1) = tr (M_2 M_1).$ The lemma is proved.

In the system (1.1) substitute
$$
\Psi = Y\Phi. \eqno (2.17)
$$
We obtain
$$
\sist{\Phi' = [A(t) + B(t) Y]\Phi,}{[Y'+ Y B(t) Y + A^*(t) Y + Y A(t) - C(t)]\Phi = 0, \ph t \ge t_0.}
$$
It follows from here and from (2.17) that all solutions of the matrix Riccati equation
$$
Y'+ Y B(t) Y + A^*(t) Y + Y A(t) - C(t) = 0, \phh t \ge t_0, \eqno (2.18)
$$
existing on any interval $[t_1,t_2) \subset [t_0,+\infty)$, are connected with solutions $(\Phi(t),\Psi(t))$ of the system (1.1) by relations.
$$
\Phi'(t) = [A(t) + B(t) Y(t)] \Phi(t), \phh \Psi(t) = Y(t) \Phi(t), \ph t \in [t_1,t_2) \eqno (2.19)
$$
(by a solution of Eq. (2.18)  on $[t_1,t_2) \subset [t_0,+\infty)$ we mean an absolutely continuous matrix function on $[t_1,t_2)$, satisfying (2.18) almost everywhere on $[t_1,t_2)$).

For any matrix $M$ of dimension $n\times n$ denote by $\lambda_1(M), \dots, \lambda_n(M)$ the eigenvalues of $M$, and if they are real-valued then we will assume that they are ordered by
$$
\lambda_1(M) \le \dots,\le \lambda_n(M).
$$
The nonnegative (positive) definiteness of any Hermitian matrix will be denoted by $H\ge~ 0  \linebreak (> 0)$. By $I$ it will be denoted the identity matrix of dimension $n\times n$.

{\bf Lemma 2.4.} {\it For a matrix $S$ of dimension $n\times n$ and any Hermitian matrix $H \ge 0$ of the same dimension the inequality
$$
tr(SHS^*) \ge \frac{\lambda_1(H)}{n}\biggl\{\biggl[tr \biggl(\frac{S + S^*}{2}\biggr)\biggr]^2 + \biggl[tr \biggl(\frac{S - S^*}{2 i}\biggr)\biggr]^2\biggr\}
$$
is valid.
}

Proof. Let $U_H$ be a unitary matrix such that
$$
U_h H U^*_H = diag \{\lambda_1(H),\dots,\Lambda_n(H)\} \stackrel{def}{=}  H_0.
$$
and let
$$
S_H = U_H S U^*_H = (s_{jk})_{j,k=1}^n.
$$
Then
$$
tr (S H S^*) = tr (S_H H_0 S_H) = \sum\limits_{j,k=1}^n \lambda_k(H) s_{jk} \overline{s_{jk}}. \eqno (2.20)
$$
Since $H\ge 0$ we have $\lambda_n(H)\ge\dots \lambda_1(H) \ge 0$. This together with (2.20) implies
$$
tr (S H S^*) \ge \lambda_1(H) \sum\limits_{j,k=1}^ns_{jk} \overline{s_{jk}} \ge \lambda_1(H)\sum\limits_{j=1}^ns_{jj} \overline{s_{jj}} =
$$
$$
 =\lambda_1(H)\biggl[\sum\limits_{j=1}^n\biggl(\frac{s_{jj} + \overline{s_{jj}}}{2}\biggr)^2 + \sum\limits_{j=1}^n\biggl(\frac{s_{jj} - \overline{s_{jj}}}{2 i}\biggr)^2\biggr] \ge
$$
$$
\ge\frac{\lambda_1(H)}{n}\biggl\{\biggl[\sum\limits_{j=1}^n\frac{s_{jj} + \overline{s_{jj}}}{2}\biggr]^2 + \sum\limits_{j=1}^n\biggl\{\biggl[\frac{s_{jj} - \overline{s_{jj}}}{2 i}\biggr]^2\biggr\} =
$$
$$
=\frac{\lambda_1(H)}{n}\biggl\{\biggl[tr \biggl(\frac{S + S^*}{2}\biggr)\biggr]^2 + \biggl[tr \biggl(\frac{S - S^*}{2 i}\biggr)\biggr]^2\biggr\}.
$$
The lemma is proved.

Let $g$ be a positive linear functional on the space of matrices of dimension $n\times n$. For any matrix $M\ge 0$ of dimension $n\times n$ set
$$
\nu_g(M)\equiv\sist{0, \ph if \ph \det M = 0,}{\{g(M^{-1})\}^{-1}, \ph if \ph \det M \ne 0.}
$$

{\bf Lemma2.5.} {\it For any matrix $M$ and any Hermitian matrix $H\ge 0$ of dimension $n\times n$ the inequality
$$
g(M^* H M) \ge \nu_g(H)[g(M)]^2 \eqno (2.21)
$$
is valid.
}

Proof. If $\det M \ne 0,$, then  $H > 0$ (since $H\ge0$),  and the inequality (2.21) is proved in [11] (see [11, p. 178]). Assume $\det M = 0$. Then since $H\ge 0$ for arbitrary small $\varepsilon > 0$ the Hermitian matrix $H_\varepsilon \equiv \varepsilon I + H$ is positive definite ($H_\varepsilon > 0$).  Therefore according to the already established fact we have
$$
g(M^* H_\varepsilon M) \ge \nu_g(H_\varepsilon) [g(M)]^2 \ge 0.
$$
Therefore
$$
g(M^* H M) = -\varepsilon g(M^* M) + g(M^* H_\varepsilon M) \ge - \varepsilon g(M^* M).
$$
From here it follows (2.21). The lemma is proved.

It is known that for every Hermitian matrix $D \ge 0$ of dimension $n\times n$ the estimates
$$
\lambda_1(D) \le g(D) \le \lambda_n(D) \eqno (2.22)
$$
are valid for every positive linear functional $g$ (see [18]). Then from the relation $g(B^{-1}(t) \ge \lambda_1(B(t)), \ph B(t) > 0$ it follows that
$$
\nu_g(B(t)) \le \lambda_1(B(t)) \le tr B(t), \phh t \ge t_0. \eqno (2.23)
$$
provided $B(t) \ge 0, \ph t \ge t_0.$ Hence, $\nu_g(B(t)), \ph t\ge t_0$ is always locally integrable for $B(t)\ge 0, \ph t \ge t_0$.

\vskip 10pt

{\bf 3. Oscillation criteria.} Hereafter by the satisfiability of a relation $\mathcal{P}$ (equality, inequality) on any interval we will mean (if it is necessary) the satisfiability of  $\mathcal{P}$ almost everywhere on that interval.

Consider the linear matrix equation
$$
B(t) X = A(t), \phh t\ge t_0. \eqno (3.1)
$$
This equation  has always  a unique solution when $B(t) > 0, \ph t \ge t_0 \ph (X = X(t) \equiv B^{-1}(t) A(t) ,\ph t \ge t_0)$. In the general case it has a solution iff (the Kronecker-Capelli theorem [12, p. 77])
$$
rank B(t) = rank (B(t) | A(t)), \phh t\ge t_0.
$$
For any matrix function $P(t), \ph t\ge t_0$ of dimension $n\times n$ set
$$
J_P(t) \equiv -\il{t_0}{t} [C(\tau) + A^*(\tau) P(\tau)]d \tau - P(t), \ph t \ge t_0.
$$
Denote by $\mathcal{M}_\mathbb{R}$ the set of matrices of dimension $n\times n$ with real entries.

{\bf Theorem 3.1.} {\it Let $A(t), \ph B(t), \ph C(t) \in \mathcal{M}_\mathbb{R}, \ph t \ge t_0$, Eq. (3.1) have a solution $F(t)$ such that $A^*(t) F(t)$ is locally integrable on $[t_0,+\infty)$, and let the following conditions be satisfied.

\noindent
I) $B(t) \ge 0, \ph t \ge t_0$.

\noindent
II) $\ilp{t_0} \nu_g(B(t)) d t = +\infty.$

\noindent
III) $\lim\limits_{t \to +\infty} g(J_F(t)) = +\infty$.

\noindent
Then the system (1.1) is oscillatory.
}

Proof. Suppose the system (1.1) is not oscillatory. Then it has a prepared solution $(\Phi(t),\Psi(t))$ such that $\det \Phi(t) \ne 0, \ph t \ge t_1$ for some $t_1 \ge t_0$.  By (2.18) and (2.19) it follows from here that for the Hermitian matrix function $Y(t) \equiv \Psi(t) \Phi^{-1}(t), \ph t\ge t_1$ the equality
$$
Y'(t) + Y(t) B(t) Y(t) + A^*(t) Y(t) + Y(t) A(t) - C(t) = 0, \phh t \ge t_1
$$
is fulfilled. Integrate this equality from $t_1$ to $t$. We obtain
$$
Y(t) - Y(t_1) + \il{t_1}{t}[Y(\tau) B(\tau) Y(\tau) + A^*(\tau) Y(\tau) + Y(\tau) A(\tau) - C(\tau)] d\tau = 0, \phh t\ge t_1. \eqno (3.2)
$$
Set $Z(t) \equiv Y(t) + F(t), \ph t \ge t_1$. Then since by the condition of the theorem $A^*(t) F(t)$ is locally integrable  from (3.2) we obtain
$$
Z(t) - Y(t_1) + \il{t_0}{t_1}[C(\tau) + A^*(\tau) F(\tau)]d\tau + \il{t_1}{t}Z^*(\tau) B(\tau) Z(\tau) d\tau + J_F(t) = 0, \eqno (3.3)
$$
$t\ge t_1$ (since by (3.1) $Z^*(\tau)B(\tau) Z(\tau) = Y(\tau) B(\tau) Y(\tau) + F^*(\tau)B(\tau) Y(\tau) + Y(\tau) B(\tau) F(\tau) + F^*(\tau) B(\tau) F(\tau) =  Y(\tau) B(\tau) Y(\tau) + A^*(\tau) Y(\tau) + Y(\tau) A(\tau) + F^*(\tau) B(\tau) F(\tau), \ph \tau \ge t_1$). By Lemma 2.5 from the condition I) it follows
$$
g(Z^*(\tau) B(\tau) Z(\tau)) \ge \nu_g(B(\tau)) [g(Z(\tau))]^2, \phh \tau \ge t_1.
$$
This together with (3.3) implies
$$
g(Z(t))-g\biggl[Y(t_1) - \il{t_0}{t_1}[C(\tau)+A^*(\tau) F(\tau)]d\tau\biggr] + \phantom{aaaaaaaaaaaaaaaaaaaaaaaaaaaaaaaaaaaaaa}
$$
$$ \phantom{aaaaaaaaaaaaaaaaaaaaaa} +\il{t_1}{t}\nu_g(B(\tau)) [g(Z(\tau))]^2d\tau+g(J_F(t)) \le 0, \phh t \ge t_1. \eqno (3.4)
$$
Without loss of generality on the basis of the condition III) we can take that $t_1$ is so large that
$$
- g\biggl[Y(t_1) - \il{t_0}{t_1}[C(\tau) + A^*(\tau) F(\tau)]d\tau\biggr] +
  g(J_F(t)) \ge 2, \ph t \ge t_1.
$$
Then from (3.4) we obtain
$$
g(Z(t)) \le -2, \phh t\ge t_1. \eqno (3.5)
$$
Set $f(t)\equiv -g(Z(t)) - \il{t_1}{t}\nu_g(B(\tau))]g(Z(\tau)]^2 d\tau, \ph f_1(t)\equiv f(t) - 1, \ph t \ge t_1.$ It follows from (3.5) that
$$
f(t) > f_1(t) > 0, \phh t \ge t_1. \eqno (3.6)
$$
Moreover $f_1(t)$ is absolutely continuous on $[t_1,+\infty)$. Consider the integral Riccati equations
$$
y(t) + \il{t_1}{t}\nu_g(B(\tau)) y^2(\tau)d\tau + f(t) =0, \ph t \ge t_1, \eqno (3.7)
$$
$$
y(t) + \il{t_1}{t}\nu_g(B(\tau)) y^2(\tau)d\tau + f_1(t) =0, \ph t \ge t_1. \eqno (3.8)
$$
It follows from (3.4) that $f_1(t) \ge -1 - g(Z(t)) -  il{t_1}{t}\nu_g(B(\tau)) [g(Z(\tau))]^2d\tau  \ge -1 - g\biggl[Y(t_1) - \il{t_0}{t_1}[C(\tau) + A^*(\tau) F(\tau)]d\tau\biggr] + g(J_F(t)), \ph t\ge t_1$. This together with the condition III) implies that
$$
\lim\limits_{t \to +\infty} f_1(t) = +\infty. \eqno (3.9)
$$
Obviously $y(t) \equiv g(Z(t)), \ph t\ge t_1$ is a solution of Eq. (3.7) on $[t_1,+\infty)$. Then in virtue of Lemma 2.2 it follows from (3.5) that Eq. (3.8) has a solution $y_1(t)$ on $[t_1,+\infty)$. Note that $y_1(t)$ is a solution of the Riccati equation
$$
y' + \nu_g(B(t)) y^2  = f_1'(t) = 0, \phh t \ge t_1.
$$
(recall that $f_1(t)$ is absolutely continuous on $[t_1,+\infty)$). Then by (2.7)  the linear system
$$
\sist{\phi' = \nu_g(B(t))\psi,}{\psi' = - f_1'(t) \phi, \ph t \ge t_1} \eqno (3.10)
$$
is not oscillatory. On the other hand since according to (3.9)
$\ilp{t_1}f_1'(\tau) d\tau - \lim\limits_{t\to +\infty}[f_1(t) - f_1(t_1)] = +\infty$ by Theorem 2.2  from the conditions I), II) it follows that the system (3.10) is oscillatory. We have obtained a contradiction, which completes the proof of the theorem.

Note that in the case when $A(t), B(t)$ and $C(t)$ are continuous and $B(t) > 0, \ph t \ge t_0$ the conditions of Theorem 3.1 become the conditions of Theorem 1.1. Therefore, Theorem~ 3.1 is a extension of Theorem 1.1. It should be noted here also that by analogy of this extension of Theorem 1.1 it can be extended Theorem 2.2 of work [11] (to do this it is needs  to substitute $Y = \alpha V, \ph \alpha > 0, \ph \alpha \in \mathbb{C}^1$ in  Eq. (2.18)).

{\bf Example 3.1.} {\it Assume $A(t) \equiv 0,\ph  B(t) = \sin^2 t  B_0, \ph C(t) \in \mathcal{M}_\mathbb{R}, \ph t \ge t_0, \linebreak \ilp{t_0}-g(C(t))d t = +\infty$, where $B_0\in \mathcal{M}_\mathbb{R}$ is a positive definite Hermitian matrix. Obviously with such $A(t), \ph B(t)$ and $C(t)$ Theorem 1.1 is not applicable to the system (1.1). For this case of $B(t)$ we have $\nu_g(B(t)) = \sin^2 t \frac{1}{g(B_0^{-1})}, \ph t\ge t_0$. Moreover $F(t) \equiv 0$ is a locally integrable solution of Eq. (3.1). Then since $\ilp{t_0}\nu_g(B(t)) = \frac{1}{g(B_0^{-1})} \ilp{t_0}\sin^2 t d t = +\infty$ and $\ilp{t_0}(-g(C(t)))d t = +\infty$ by Theorem 3.1 the system (1.1) for the considered case of its coefficients is oscillatory.
}

{\bf Theorem 3.2.} {\it Let $A(t), \ph B(t), \ph C(t) \in \mathcal{M}_\mathbb{R}, \ph B(t) \ge 0, \ph t \ge t_0$, and let  $F(t)$ be an absolutely continuous solution of Eq. (3.1). If for some positive linear functional $g$ the scalar system
$$
\sist{\phi. =\nu_g(B(t)) \psi,}{\psi' =-g[C(t) + A^*(t) F(t) + F'(t)] \phi, \ph t \ge t_0} \eqno (3.11)
$$
is oscillatory, then the system (1.1) is also oscillatory.
}

Proof. Suppose the system (1.1) is not oscillatory. Then there exists its a prepared solution $(\Phi(t),\Psi(t))$
such that $\det \Phi(t) \ne 0, \ph t \ge t_1$, for some $t_1 \ge t_0$. By (2.18) and (2.19) it follows from here that for the Hermitian matrix function  $Y(t) \equiv \Psi(t) \Phi^{-1}(t), \ph t \ge t_1$ the equality
$$
Y'(t) + Y(t) B(t) Y(t) + A^*(t) Y(t) + Y(t) A(t) - C(t) = 0,, \phh t \ge t_1
$$
is fulfilled. If we set $Z(t) = Y(t) + F(t), \ph t \ge t_1$, then from the above equality we obtain
$$
Z'(t) + Z^*(t) B(t) Z(t) - C(t) - A^*(t) F(t) - F'(t) = 0,  \phh t \ge t_1. \eqno (3.12)
$$
Since $B(t) \ge 0, \ph t \ge t_0$ by virtue of Lemma 2.5 we have
$$
g[Z^*(t) B(t) Z(t)] \ge \nu_g(t) [g(Z(t))]^2, \ph t \ge t_1.
$$
This together with (3.12) implies
$$
[g(Z(t)]' + \nu_g(B(t))[g(t)]^2 - g[C(t) + A^*(t) F(t) + F'(t)] \le 0, \phh t \ge t_1.  \eqno (3.13)
$$
Set $f(t) \equiv - [g(Z(t))]' - \nu_g(B(t))[g(t)]^2, \ph t \ge t_1$. It follows from (3.13) that
$$
f(t) \ge - g[C(t) + A^*(t) F(t) + F'(t)], \phh t \ge t_1.  \eqno (3.14)
$$
Consider the scalar Riccati equations
$$
y' + \nu_g(B(t)) y^2 - g[C(t) + A^*(t) F(t) + F'(t)] = 0, \ph t \ge t_1, \eqno (3.15)
$$
$$
y' + \nu_g(B(t)) y^2 + f(t) = 0, \phh t \ge t_1. \eqno (3.16)
$$
Obviously $y(t) \equiv g(Z(t)), \ph t \ge t_1$ is a solution  of Eq. (3.16) on $[t_1,+\infty)$. Then applying Theorem 2.1 to the pair of equations (3.15) and (3.16), and taking into account (3.14) we conclude that Eq. (3.15) has a solution on $[t_1,+\infty)$. By (2.7) it follows from here that the system (3.11) is not oscillatory, which contradicts the condition of the theorem. The obtained contradiction completes the proof of the theorem

Note that if under the restriction, that $F(t)$ is absolutely continuous on $[t_0,+\infty)$, the conditions I) - III) are satisfied, then by Theorem 2.2 the system (3.11) is oscillatory. Hence, Theorem 3.2 is a complement to Theorem 3.1 (therefore to Theorem 1.1).

Denote
$$
J(t) \equiv tr \biggl[\frac{A(t) + A^*(t)}{2}B^{-1}(t)\biggr] - \il{t_0}{t}tr \Bigl[A(\tau) B^{-1}(\tau) A^*(\tau) + C(\tau)\Bigr] d \tau +
 $$
 $$
 +\il{t_0}{t}\frac{\lambda_1(B(\tau))}{n}\biggl[ tr\biggl(\frac{A(\tau) - A^*(\tau)}{2 i}\biggr)\biggr]^2 d \tau, \ph t \ge t_0.
$$

{\bf Theorem 3.3.} {\it Let the the functions  $tr \Bigl[A(t) B^{-1}(t) A^*(t)\Bigr], \ph \lambda_1(B(t))\biggl[ tr \bigl(A(t) - A^*(t)\bigr)\biggr]^2, \linebreak   t \ge t_0$ be locally integrable and let the following conditions be satisfied.

\noindent
I') $B(t) > 0, \ph t \ge t_0$.

\noindent
IV) $\ilp{t_0} \lambda_1(B(\tau)) d \tau = \lim\limits_{t \to +\infty} J(t) = +\infty.$

\noindent
Then the system (1.1) is oscillatory.
}

Proof. Suppose the system (1.1) is not oscillatory. Then by (2.19) Eq. (2.18) has a Hermitian solution $Y(t)$ on $[t_1,+\infty)$ for some $t_1 \ge t_0$. Then using I') we can write
$$
Y'(t) + \frac{1}{2}\Bigl\{[Y(t) + A(t) B^{-1}(t)] B(t) [Y(t) + B^{-1}(t) A^(t)] +
$$
$$
 +[Y(t) + A^*(t) B^{-1}(t)] B(t) [Y(t) + B^{-1}(t) A(t)]\Bigr\}+
$$
$$
+\frac{A^*(t) - A(t)}{2} Y(t) + Y(t) \frac{A(t) - A^*(t)}{2} -
$$
$$
-\frac{1}{2}[A(t)B^{-1}(t) A^*(t) + A^*(t) B^{-1}(t) A(t)] - C(t) = 0, \ph t\ge t_1. \eqno (3.17)
$$
Since $Y(t)$ and $B^{-1}(t)$ are Hermitian by Lemma 2.4 we have
$$
tr \frac{1}{2}\Bigl\{[Y(t) + A(t) B^{-1}(t)] B(t) [Y(t) + B^{-1}(t) A^(t)] +
$$
$$
 +[Y(t) + A^*(t) B^{-1}(t)] B(t) [Y(t) + B^{-1}(t) A(t)]\Bigr\} \ge
$$
$$
\frac{\lambda_1(B(t))}{2n}\biggl\{\biggl[tr\biggl(Y(t) + \frac{A(t) B^{-1}(t) + B^{-1}(t) A^*(t)}{2}\biggr)\biggr]^2 + \biggl[tr\biggl(\frac{A(t) B^{-1}(t) - B^{-1}(t) A^*(t)}{2 i}\biggr)\biggr]^2
$$
$$
+\biggl[tr\biggl(Y(t) + \frac{A^*(t) B^{-1}(t) + B^{-1}(t) A(t)}{2}\biggr)\biggr]^2 + \biggl[tr\biggl(\frac{A^*(t) B^{-1}(t) - B^{-1}(t) A(t)}{2 i}\biggr)\biggr]^2\biggr\},
$$
$t \ge t_1$. By Lemma 2.3 from here we obtain
$$
tr \frac{1}{2}\Bigl\{[Y(t) + A(t) B^{-1}(t)] B(t) [Y(t) + B^{-1}(t) A^(t)] +
$$
$$
 +[Y(t) + A^*(t) B^{-1}(t)] B(t) [Y(t) + B^{-1}(t) A(t)]\Bigr\} \ge
$$
$$
\frac{\lambda_1(B(t))}{n}\biggl\{\biggl[tr\biggl(Y(t) + \frac{A(t) +  A^*(t)}{2}  B^{-1}(t)\biggr)\biggr]^2 + \biggl[tr\biggl(\frac{A(t) - A^*(t)}{2i} B^{-1}(t)\biggr)\biggr]^2\biggr\}, \ph t \ge t_1.
$$
This together with (3.17) implies
$$
tr Y'(t) + \frac{\lambda_1(B(t))}{n}\biggl\{\biggl[tr\biggl(Y(t) + \frac{A(t) +  A^*(t)}{2}  B^{-1}(t)\biggr)\biggr]^2 - \phantom{aaaaaaaaaaaaaaaaaaaaaaaaaaaaa}
$$
$$
-\frac{1}{2}tr\Bigl\{[A(t)B^{-1}(t) A^*(t) + A^*(t) B^{-1}(t) A(t)] - C(t)\Bigr\} + \biggl[tr\biggl(\frac{A(t) - A^*(t)}{2i} B^{-1}(t)\biggr)\biggr]^2\biggr\} \le 0,
$$
$t \ge t_1$.  If we substitute $Z(t) \equiv Y(t) + \frac{A(t) + A^*(t)}{2} B^{-1}(t), \ph t \ge t_1$ in the above inequality  and integrate (by taking into account the condition of local integrability of the functions $tr [A(t)B^{-1}(t) A^*(t)], \ph tr \lambda_1(B(t))[tr(A(t) - A^*(t))]^2$) from $t_1$ to $t$ we obtain
$$
tr Z(t) + \il{t_1}{t}\frac{\lambda_1(B(\tau))}{n}[tr Z(\tau)]^2 d \tau + J(t) + c \le 0, \phh t \ge t_1,
$$
where $c = Y(t_1) + \il{t_0}{t_1} tr [C(\tau) + A(\tau) B^{-1}(\tau) A^*(\tau)] d \tau- \il{t_0}{t_1}\biggl[tr\biggl(\frac{A(\tau) - A^*(\tau)}{2i} B^{-1}(\tau)\biggr)\biggr]^2d\tau = const.$ Further as in the proof of Theorem 3.1. The theorem is proved.

Set
$$
\nu_0(B(t)) \equiv \sist{0, \ph if \ph \det B(t) = 0,}{\frac{1}{tr (B^{-1}(t))}, \ph if \ph \det B(t) \ne 0,} \ph t \ge t_0.
$$
It is not difficult to verify that
$$
\frac{1}{tr(B^{-1}(t))} \le \lambda_1(B(t)) \le \frac{n}{tr(B^{-1}(t))}, \ph \mbox{for all} t \ge t_0, \ph \mbox{for whicch} \ph B(t) > 0. \eqno (3.18)
$$
Therefore Theorem 3.3 remains valid if we replace the condition $\ilp{t_0}\lambda_1(B(t)) d t = +\infty$ of Theorem 3.3 by the following one $\ilp{t_0}\nu_0(B(t)) d t = +\infty$. Moreover by (2.22) in the case
$-\il{t_0}{t} \Bigl(C(s) + A^*(s) B^{-1}(s) A(s)\Bigr) ds - B^{-1}(t) A(t) \ge 0, \ph t \ge T$, for some $T\ge t_0$,
the functional $g$ in Theorem 1.1 is equivalent to the functional $tr$. Hence due to (3.18) Theorem 3.3 is a complement to Theorem 1.1.

{\bf Theorem 3.4.} {\it Let  the following conditions be satisfied.

\noindent
I') $B(t) > 0, \ph t \ge t_0.$

\noindent
V) $\ilp{t_0}\frac{d t}{tr (B^{-1}(t))} = +\infty.$

\noindent
VI)  the  function   $tr [(A(t) + A^*(t)) B^{-1}(t)(A(t) + A^*(t))], \ph t \ge t_0$ is locally integrable on $[t_0,+\infty)$ and  $\lim\limits_{t \to +\infty}- tr \biggl[2(A(t) + A^*(t)) B^{-1}(t) +$

\phantom{aaaaaaaaaaaaaaaaa} $+\il{t_0}{t}\Bigl((A(\tau) + A^*(\tau)) B^{-1}(\tau) (A(\tau) + A^*(\tau)) + 4 G(\tau)\Bigr) d\tau\Biggr] = +\infty$.

\noindent
Then the system (1.1) is oscillatory.
}

Proof. Suppose the system (1.1) is not oscillatory. Then by (2.19) Eq. (2.18) has a solution $Y(t)$ on $[t_1,+\infty)$ for some $t_1 \ge t_0$. Hence,
$$
Y'(t) + Y(t) B(t) Y(t) + A^*(t) Y(t) + Y(t) A(t) - C(t) = 0, \phh t \ge t_1.
$$
From here it follows
$$
tr \biggl\{Y'(t) + \Bigl[Y(t) + \frac{A(t) + A^*(t)}{2} B^{-1}(t)\Bigr] B(t) \Bigl[Y(t) +  B^{-1}(t) \frac{A(t) + A^*(t)}{2}\Bigr] +
$$
$$
+\frac{A^*(t) - A(t)}{2} Y(t) + Y(t) \frac{A(t) - A^*(t)}{2} - \frac{A^*(t) + A(t)}{2} B^{-1}(t)\frac{A^*(t) + A(t)}{2} - C(t)\biggr\} = 0,
$$
$t \ge t_1$.  Substitute $Z(t) \equiv Y(t)  + \frac{A(t) + A^*(t)}{2}, \ph t \ge t_1$ in the obtained equality and integrate from $t_1$ to $t$. Taking into account the fact that the function $tr [(A(t) + A^*(t)) B^{-1}(t)(A(t) + A^*(t))]$ is locally integrable, we obtain
$$
tr \biggl\{Z(t)  + \il{t_1}{t}Z(\tau) B(\tau) Z^*(\tau) d \tau + \il{t_1}{t}\Bigl[\frac{A^*(\tau) - A(\tau)}{2} Y(\tau) + \phantom{aaaaaaaaaaaaaaaaaaaaaaaaaaaaaaaaaa}
$$
$$
\phantom{aaaaaaaaaaaaaaaaaaaaaa}+Y(\tau)\frac{A(\tau) - A^*(\tau)}{2}\Bigr] d\tau + J_1(t)\biggr\} = 0, \phh t \ge t_1, \eqno (3.19)
$$
 where
$$
J_1(t) \equiv -Y(t_1) - \frac{A(t) + A^*(t)}{2} B^{-1}(t) - \phantom{aaaaaaaaaaaaaaaaaaaaaaaaaaaaaaaaaaaaaaaaaaaaaa}
$$
$$
\phantom{aaaaaaaaaaaaaaaaaaaaaaaa}-\il{t_1}{t}\Bigl[\frac{A(\tau) + A^*(\tau)}{2} B^{-1}(\tau)\frac{A(\tau) + A^*(\tau)}{2} + C(\tau)\Bigr] d\tau, \ph t \ge t_1.
$$
By Lemma 2.3 we have
$$
tr \biggl[\il{t_1}{t}\Bigl[\frac{A^*(\tau) - A(\tau)}{2} Y(\tau) + Y(\tau)\frac{A(\tau) - A^*(\tau)}{2}\Bigr] d\tau\biggr] = 0, \phh t \ge t_1. \eqno (3.20)
$$
Since $Y(t) , \ph A(t) + A^*(t)$ and $B^{-1}(t)$ are Hermitian we have also $tr (Z(t) - Z^*(t)) = 0, \linebreak t \ge~ t_1$. By Lemma 2.4 from here we obtain
$$
tr \il{t_1}{t} Z(\tau) B(\tau) Z^*(\tau) d \tau  \ge \il{t_1}{t}\frac{\lambda_1(B(\tau))}{n}\Bigl[ tr \frac{Z(\tau) + Z^*(\tau)}{2}\Bigr]^2 d \tau  = \il{t_1}{t}\frac{\lambda_1(B(\tau))}{n}\Bigl[ tr Z(\tau)\Bigr]^2 d \tau,
$$
$t \ge t_1$. This together with (3.19) and (3.20) implies that
$$
tr Z(t) - tr Z(t_1) + \il{t_1}{t}\frac{\lambda_1(B(\tau))}{n}\Bigl[ tr Z(\tau)\Bigr]^2 d \tau + tr J_1(t) \le 0, \ph t \ge t_1.
$$
Further as in the proof of Theorem 3.1 one can show that, if the conditions I'), VI) and the condition

\noindent
V') $\ilp{t_0}\lambda_1(B(t)) d t = +\infty$

\noindent
are satisfied, then the system (1.1) is oscillatory. But according to (3.18) the condition V') with I') is equivalent to the condition V). Therefore under the conditions of the theorem the system (1.1) is oscillatory. The theorem is proved.

{\bf Example 3.2.} {\it Let $B(t) \equiv I, \ph A(t) \equiv A_0, \ph C(t) \equiv - A^*_0 A_0,  \ph t \ge t_0,\ph
 A_0= const$ is a real-valued matrix of dimension $n\times n$.
Then
$$
\lim\limits_{t \to +\infty} g\biggl[-\il{t_0}{t}\biggl(C(\tau) + A^*(\tau) B^{-1}(\tau) A(\tau)\biggr) d \tau - B^{-1}(t) A(t)\biggr] =\lim\limits_{t \to+\infty}g[-A_0] \ne +\infty.
$$
Therefore, for this particular case Theorem 1.1 is not applicable to the system (1.1).  It is not difficult to verify that
$$
tr C(t) = - tr (A_0^* A_0) < 0, \ph t \ge t_0.
$$
Then since $A_0 + A^*_0 = 0$ using Theorem 3.4 to the system (1.1) we conclude that for this particular case the system (1.1) is oscillatory.}

Assume $B(t) \ge 0, \ph t \ge t_0$ and let $\sqrt{B(t)}, \ph t \ge t_0$ be absolutely continuous. Consider the linear matrix equation
$$
\sqrt{B(t)} X (A(t)\sqrt{B(t)} - \sqrt{B(t)}') = A(t)\sqrt{B(t)} - \sqrt{B(t)}', \ph t \ge t_0. \eqno (3.21)
$$
This equation has always a solution when $B(t) > 0, \ph t \ge t_0$. But it can have also a solution when $B(t)$ is not invertible for all (for some) $t \ge t_0$ (see [5]). In the general case Eq. (2.21) has a solution if and only if the equations
$$
\sqrt{B(t)}Y = A(t)\sqrt{B(t)} - \sqrt{B(t)}', \ph t \ge t_0, \phantom{aaaaaaaaaaaaaaaaaa}
$$
$$
 \phantom{aaaaaaaaaaaaaaaaaa} Z  A(t)\sqrt{B(t)} - \sqrt{B(t)}' =  A(t)\sqrt{B(t)} - \sqrt{B(t)}', \ph t\ge t_0
$$
have solutions (see [4], p. 23). Hence, Eq. (3.21) has a solution if and only if \linebreak $rank \sqrt{B(t)} = rank (\sqrt{B(t)}|  A(t)\sqrt{B(t)} - \sqrt{B(t)}'), \ph t \ge t_0.$

Let $F(t)$ be a solution of Eq. (3.21). We set:
$$
A_F(t) \equiv F(t)(A(t) \sqrt{B(t)} - \sqrt{B(t)}'), \phh J_2(t)\equiv -\frac{1}{2}tr (A_F(t) + A^*_F(t)) -
$$
$$
-\il{t_0}{t} tr \biggl[A_F(\tau) A^*_F(\tau)   + B(\tau) C(\tau)\biggr]d \tau +
\frac{1}{n}\il{t_0}{t}\Bigl[tr\Bigl(\frac{ A_F(\tau) - A^*_F(\tau)}{2i}\Bigr)\Bigr]^2 d\tau, \phh t \ge t_0.
$$

{\bf Theorem 3.5.} {\it Let $\sqrt{B(t)}$ be absolutely continuous on $[t_0,+\infty)$ and let $F(t)$ be a solution of Eq. (3.10) such that the functions  $tr (A_F(t) A^*_F(t) + B(t) C(t)), \ph [tr( A_F(t) - A^*_F(t))]^2, \ph t \ge t_0$ are locally integrable on $[t_0,+\infty)$. If
$$
\lim\limits_{t\to+\infty} J_2(t) = +\infty
$$
then the system (1.1) is oscillatory.
}

Proof. Suppose the system (1.1) is not oscillatory. Then by (2.19) Eq. (2.18) has a solution $Y(t)$ on $[t_1,+\infty)$ for some $t_1 \ge t_0$. Hence,
$$
Y'(t) + Y(t)B(t)Y(t) + A^*(t) Y(t) + Y(t) A(t) - C(t) = 0, \phh t \ge t_1.
$$
Multiply both sides of this equality at left and at right by $\sqrt{B(t)}$. Taking into account the equality
$$
(\sqrt{B(t)} Y(t) \sqrt{B(t)})' = \sqrt{B(t)}' Y(t) \sqrt{B(t)} + \sqrt{B(t)} Y'(t) \sqrt{B(t)} + \sqrt{B(t)} Y(t) \sqrt{B(t)}',
$$
$t \ge t_1$ we obtain
$$
(\sqrt{B(t)} Y(t) \sqrt{B(t)})' + (\sqrt{B(t)} Y(t) \sqrt{B(t)})^2 + (\sqrt{B(t)} A^*(t) - \sqrt{B(t)}') Y(t) \sqrt{B(t)} +\phantom{aaaaaaa}
$$
$$
\phantom{aaaaaaaaaa}+ \sqrt{B(t)} Y(t) (A(t)\sqrt{B(t)} - \sqrt{B(t)}') - \sqrt{B(t)} C(t)\sqrt{B(t)}, \phh t \ge t_1.  \eqno (3.22)
$$
Since $F(t)$ is a solution of Eq. (3.21) we have $\sqrt{B(t)} A^*(t) - \sqrt{B(t)}' = (\sqrt{B(t)} A^*(t) - \sqrt{B(t)}') F^*(t) \sqrt{B(t)} = A^*_F(t)\sqrt{B(t)}= A^*_F(t) \sqrt{B(t)}, \ph A(t)\sqrt{B(t)} - \sqrt{B(t)}' = A_F(t) \sqrt{B(t)},\linebreak  t \ge t_1$. From here and from (3.22) it follows
$$
tr \{V(t)\}' + tr \{V^2(t) + A^*_F(t) V(t) + V(t) A_F(t) - \sqrt{B(t)} C(t) \sqrt{B(t)}\} =0, \ph t \ge t_1, \eqno (3.23)
$$
 where $V(t) \equiv \sqrt{B(t)} Y(t) \sqrt{B(t)}, \ph t \ge t_1.$
By Lemma 2.3
$$
tr[\sqrt{B(t)} C(t) \sqrt{B(t)}] = tr[B(t) C(t)], \phh t \ge t_0.
$$
Then if we substitute $V(t) \equiv Z(t) - \frac{A_F(t) + A^*_F(t)}{2}, \ph t \ge t_1$ in (3.23) (except in the expression $tr\{V(t)\}'$)  and take into account the condition of local integrability of $tr [A_f(t) A^*_F(t) + B(t) C(t)]$ and $[tr (A_F(t) - A^*_F(t))]^2$ we can, as in the proof of Theorem 3.3, to derive the inequality
$$
[tr Z(t)] + \il{t_1}{t}[tr Z(\tau)]^2 + J_2(t) + c_1 \le 0,
$$
where $c_1 \equiv - tr V(t_1) + \il{t_0}{t_1} tr [A_F(t) A^*_F(t) + B(t) C(t)] - \frac{1}{4 n} \il{t_0}{t_1} [tr (A_F(t) - A^*_F(t))]^2$  is a constant. Further as in the proof of Theorem 3.1. The theorem is proved.

{\bf Example 3.3.} {\it Assume $B(t) \equiv \begin{pmatrix} I_m & \theta_{12}\\ \theta_{21} & \theta_{22}\end{pmatrix}, \ph A(t)\equiv \begin{pmatrix} \theta_{11} & A_1(t)\\ \theta_{21} & A_2(t)\end{pmatrix}, \ph C(t)\equiv I, \linebreak rank A_2(t) \not\equiv 0, \ph t \ge t_0,$  where $I_m$ is an identity matrices of dimensions $m\times m$ ($m < n$), $\theta_{11}, \ph \theta_{12}, \ph \theta_{21}$ and $\theta_{22}$ are null matrices of dimensions $m\times m, \ph (n-m)\times m, \ph m\times (n-m)$ and $(m-m)\times (n-m)$ respectively. Obviously $rank B(t)\not\equiv rank (B(t) | A(t)), \ph t \ge t_0.$ Therefore Eq. (3.1) has no solution, which means that for this particular case Theorems 3.1 and 3.2 are not applicable to the system (1.1). Obviously for this case $F(t)\equiv 0, \ph t \ge t_0$ is a solution for Eq. (3.21). Then $A_F(t)\equiv 0, \phh t \ge t_0,, \ph J_2(t) =  (t - t_0) m \to +\infty$ for $t \to +\infty.$. By Theorem 3.5 it follows from here that in this particular case the system (1.1) is oscillatory.
}

Denote by $\Omega_n$ the set of $n\times n$ dimensional matrices $M$ for which
$$
Re \hskip 2pt \lambda_1(M) =\dots =Re \hskip 2pt \lambda_n(M).
$$
Let $\Lambda(t) \in \Omega_n, \ph t \ge t_0$ be a complex-valued locally integrable matrix function on $[t_0,+\infty)$. Consider the linear matrix equation
$$
B(t) X + X B(t) = \Lambda(t) + \Lambda^*(t) + A(t) + A^*(t), \phh t\ge t_0. \eqno (3.24)
$$
Note that if $X(t), \ph t \ge t_0$ is any solution of this equation, then $H(t)\equiv\frac{X(t) + X^*(t)}{2}, \ph t \ge t_0$ is its a Hermitian solution. Indicate some particular cases, when Eq. (3.24) has a solution.

\noindent
I$^\circ$) $B(t) > 0, \ph t \ge t_0$. In this case $B(t)$ and $-B(t)$ have no common eigenvalues. Then (see [4], pp. 203, 207)
Eq. (3.24) has a unique (therefore Hermitian) solution, which can be given in the following closed form (see[2], p. 212, Theorem 6)
$$
H_\Lambda(t) \equiv \ilp{0}\exp\biggl\{-\tau B(t)\biggr\}\biggl[\Lambda(t) +\Lambda^*(t) + A(t) + A^*(t)\biggr]\exp\biggl\{-\tau B(t)\biggr\} d \tau, \phh t \ge t_0.
$$
Note that this integral converges and gives a hermitian solution for Eq. (3.24) not only  for the case $B(t) > 0, \ph t \ge t_0$, but also for a more general case, when $B(t) \ge 0, \ph t \ge t_0$ and $\Lambda(t) + \Lambda^*(t) + A(t) + A^*(t) = 0$ for all $t \ge t_0$ for which $\lambda_1(B(t)) = 0$.

\noindent
II$^\circ$) $rank B(t) \ge n-1, \ph t \ge t_0.$ Show that  there exists a real-valued  locally integrable function $\mu(t), \ph t \ge t_0$ such that if for some $\Lambda(t) \in \Omega_n, \ph t\ge t_0 \ph \Lambda(t) + \Lambda^*(t) = \mu(t)I, \ph t \ge t_0,$ then Eq. (3.24) has a solution. If $rank B(t) = n$, then we have the considered case I). Suppose $rank B(t) = n-1 \ph  (t$ is fixed). Let $U (t)$ be a $n\times n$ dimensional unitary matrix  such that
$$
U(t) B(t) U^*(t) = \ diag\{b_1(t),\dots,b_n(t)\}\stackrel{def}{=} B_0(t), \phh 0 = b_1(t) < b_2(t)\le\dots\le b_n(t).
$$
Then Eq. (3.24) is equivalent to the following
$$
B_0(t) V + V B_0(t) = \mathcal{A}(t), \eqno (3.25)
$$
where $V\equiv U(t) X U^*(t), \phh \mathcal{A}(t)\equiv U(t)[\Lambda(t) + \Lambda^*(t) + A(t) + A^*(t)] U^*(t).$ If we write
$$
B_0(t)=\begin{pmatrix}0 & \theta\\ \theta^T & B_1(t)\end{pmatrix}, \ph V = \begin{pmatrix}0 & v_{12}\\ v_{21} & V_{22}\end{pmatrix}, \ph \mathcal{A}(t) = \begin{pmatrix}a_{11}(t) & a_{12}(t)\\ a_{21}(t) & A_{22}(t)\end{pmatrix}, \ph t \ge t_0,
$$
where $\theta \equiv (0,\dots,0)$ is a null vector of dimension $n-1, \ph \theta^T$ is the transpose of $\theta, \ph v_{12}$ and $a_{12}(t)$ are matrices of dimension $1\times n, \ph v_{21}$ and $a_{2`1}(t)$ are matrices of dimension $n\times 1, V_{22}$ and $A_{22}(t)$ are matrices of dimension $(n-1) \times (n-1)$, then Eq. (3.25) can be vritten equivalently in the form
$$
\left\{\begin{array}{l}
v_{12}B_1(t) = a_{12}(t),\\
B_1(t) v_{21}(t) = a_{21}(t),\\
B_1(t)V_{22} + V_{22}B_1(t) = A_{22}(t),\\
2 Re \hskip 2pt \mu(t) + a_{11}(t) = 0.
\end{array}
\right. \eqno (3.26)
$$
Since $rank B(t) = n-1$ we have $B_1(t) > 0$. Moreover since $A(t)$ is Hermitian $A_{11}(t)$ is real valued. Then, obviously, for $\mu(t)\equiv -\frac{1}{2} a_{11}(t), \ph t \ge t_0$ the system (3.26) has a solution. Thus Eq. (3.24) has always a solution with this $\mu(t)$ and the appropriate chosen $\Lambda(t)$.

\noindent
III$^\circ$) $B(t)\equiv \begin{pmatrix}\theta_{11} & \theta_{12}\\ \theta_{21}& B_{22}(t)\end{pmatrix}, \ph \Lambda(t) + A(t) \equiv\begin{pmatrix}\theta_{11} & \theta_{12}\\ A_{21}(t) & A_{22}(t)\end{pmatrix}, \ph t \ge t_0,$ where $\theta_{11}, \ph \theta_{12}$ and $\theta_{21}$ are null matrices of dimensions $m\times m, \ph (n-m)\times m$ and $(n-m)\times (n-m)$ respectively, $A_{21}(t)$ is a matrix function of dimension $m\times (n-m), \ph B_{22}(t)$ and $A_{22}(t)$ are matrix functions of dimension $(n-m)\times (n-m), \ph (0 < m < n)$. If we write $X\equiv \begin{pmatrix}x_{11} & x_{12}\\x_{21} & x{22}\end{pmatrix}$, then the system (3.24) we can rewrite in the form of the following system
$$
\left\{\begin{array}{l}
x_{11} = 0,\\
x_{12}B_{22}(t) A^*_{21}(t),\\
B_{22}(t) x_{21} = A_{21}(t),\\
B_{22}(t) x_{22} + x_{22} B_{22}(t) = \Lambda(t) + \Lambda^*(t) + A(t) + A^*(t),
\end{array}
\right.  \phh t \ge t_0.
$$
Then if $B_{22}(t) > 0, \ph t \ge t_0$, the Hermitian solution $H(t)$ of Eq. (3.24) can be given in the closed form
$H(t) =\begin{pmatrix}0 & H_{12}(t)\\H_{12}^*(t) & H_{22}(t)\end{pmatrix}, \ph t \ge t_0$, where $H_{12}(t)\equiv A_{21}^*(t) B_{22}^{-1}(t),$
$$
H_{22}(t) \equiv\ilp{0}\exp\biggl\{-\tau B_{22}(t)\biggr\}\biggl[\Lambda(t) + \Lambda^*(t) + A(t) + A^*(t)\biggr]\exp\biggl\{-\tau B_{22}(t)\biggr\} d \tau, \phh t \ge t_0.
$$

For any absolutely continuous matrix function $F(t)$ of dimension $n\times n$ on $[t_0,+\infty)$ set
$$
\mathcal{D}_F(t) \equiv - F'(t) + F(t) B(t) F(t) - F(t) A(t) - A^*(t) F(t) - C(t), \phh t \ge t_0.
$$

{\bf Theorem 3.6.} {\it Let the following conditions be satisfied.

\noindent
I) $B(t) \ge 0, \ph t \ge t_0$

\noindent
VII) For some locally integrable $\Lambda(t) \in \Omega_n, \ph t \ge t_0$ Eq. (3.24) has an solution $F(t)$ such that $\mathcal{D}_F(t)$ is locally integrable on $[t_0,+\infty)$.

\noindent
VIII) The scalar system
$$
\sist{\frac{1}{n} tr [\Lambda(t) + \Lambda^*(t)] \phi + \frac{\lambda_1(B(t))}{n} \psi,}{\psi' = -[tr \mathcal{D}_F(t)] \phi, \ph t \ge t_0}
$$
is oscillatory.

\noindent
Then the system (1.1) is oscillatory.
}

Proof. Suppose the system (1.1) is not oscillatory. Then it has a conjoined solution $(\Phi(t),\Psi(t))$ such that $\det \Phi(t) \ne 0, \ph t \ge t_1$ for some $t_1\ge t_0$. By (2.19) form here it follows  that the Hermitian matrix function $Y(t)\equiv \Psi(t) \Phi^{-1}(t), \phh t \ge t_1$ is a solution for Eq. (2.18) on $[t_1,+\infty)$. Then it is not difficult to verify that for the matrix function $Z(t) \equiv Y(t) + F(t), \ph t \ge t_1$ the equality
$$
Z'(t) + Z(t) B(t) Z(t) + (A^*(t) - F(t)B(t)) Z(t) + Z(t) (A(t) - B(t) F(t)) + \mathcal{D}_F(t) = 0, \eqno (3.27)
$$
$t\ge t_1$.  Since $tr$ is a positive linear functional we have
$$
tr [Z(t) B(t) Z(t)] \ge \frac{\lambda_1(B(t))}{n}[tr Z(t)]^2, \phh t \ge t_1. \eqno (3.28)
$$
By virtue of Lemma 2.3 we have
$$
tr[(A^*(t) - F(t)B(t)) Z(t) + Z(t) (A(t) - B(t) F(t))] = \phantom{aaaaaaaaaaaaaaaaaaaaaaaaaaaaaaaaaa}
 $$
 $$
\phantom{aaaaaaaaaaaaaaaaaaaaaaaaa} +tr [Z(t) (A(t) + A^*(t) - B(t) F(t) - F(t) B(t))], \phh t \ge t_1.
$$
Then by (3.24) from here we obtain
$$
tr[(A^*(t) - F(t)B(t)) Z(t) + Z(t) (A(t) - B(t) F(t))] = tr [Z(t)(\Lambda(t) + \Lambda^*(t))], \ph t \ge t_1. \eqno (3.29)
$$
Since $\Lambda(t) \in \Omega_n, \ph t \ge t_0$ it is not difficult to verify that
$$
tr [Z(t)(\Lambda(t) + \Lambda^*(t))] = \frac{1}{n} tr [\Lambda (t) + \Lambda^*(t)] tr Z(t), \phh t \ge t_1.
$$
This together with (3.27)-(3.29) implies that
$$
[tr Z(t)]' + \frac{\lambda_1(B(t))}{n}[tr Z(t)]^2 + \frac{1}{n} tr[\Lambda(t) + \Lambda^*(t)] tr Z(t) + tr \mathcal{D}_F(t) \le 0, \phh t \ge t_1.
$$
Further as in the proof of Theorem 3.2. The theorem is proved.

{\bf Remark 3.1.} {\it Using (2.23) on the basis of comparison Theorem 2.1 one can show that in Theorem 3.6 $\lambda_1(B(t))$ can be replaced by a more explicit function $\nu_0(B(t))$.
}

If for some locally integrable matrix function $\Lambda(t) \in \Omega_n, \ph t \ge t_0$ the matrix function $\Lambda(t) + A(t), \ph t \ge t_0$ is skew symmetric ($\Lambda(t) + A(t) = - \Lambda^*(t) - A^*(t), \ph t \ge t_0$), then $F(t)\equiv 0, \ph t \ge t_0$ is a Hermitian solution of Eq. (3.24) and, hence, $\mathcal{D}_F(t) = - C(t), \ph t \ge t_0$. Therefore combining Theorem 3.6 with Theorem 2.2  and taking into account Remark 3.1 we obtain immediately

{\bf Corollary 3.1.} {\it Let the following conditions be satisfied

\noindent
I) $B(t) \ge 0, \ph t \ge t_0$.

\noindent
IX) There exists a locally integrable on $[t_0,+\infty)$ matrix function $\Lambda(t) \in \Omega_n, \ph t \ge t_0$ such that $\Lambda(t) + A(t), \ph t \ge t_0$ is skew symmetric.

\noindent
$$
X) \ph \ilp{t_0}\nu_0(B(t))\exp\biggl\{-\frac{1}{n}\il{t_0}{t}tr[\Lambda(\tau) + \Lambda^*(\tau)]d\tau\biggr\} d t = \phantom{aaaaaaaaaaaaaaaaaaaaaaaaaaaaaaa}
$$
$$
\phantom{aaaaaaaaaaaaaaaaaaaaaaaaaa} = -\ilp{t_0}tr [C(t)]\exp\biggl\{\frac{1}{n}\il{t_0}{t}tr[\Lambda(\tau) + \Lambda^*(\tau)]d\tau\biggr\} d t = +\infty.
$$

\noindent
Then the system (1.1) is oscillatory.
}

$\phantom{aaaaaaaaaaaaaaaaaaaaaaaaaaaaaaaaaaaaaaaaaaaaaaaaaaaaaaaaaaaaaaaaaaaa} \blacksquare$

\vskip 5pt

{\bf Remark 3.2} {\it Corollary 3.1 is a generalization of Lighton's oscillation criterion (see [15, Theorem 2.24]).
}

For any matrix $L\equiv (l_{jk})_{j,k =1}^n$ denote $Sum(L)\equiv \sum\limits_{j,k=1}^n l_{jk}$.
Define the Hermitian matrix $H = H_L = (h_{jk})_{j,k =1}^n$ by elements of matrix  $L=(l_{jk})_{jk=1}^n$ by formulae:

\noindent
$h_{kk} = - \sum\limits_{j=1}^n Re\hskip 2pt l_{jk}, \ph k=\overline{1,n};$

\noindent
$h_{nk}= \overline{h_{kn}} = - i \sum\limits_{j=1}^n Im \hskip 2pt l_{jk} + \frac{i}{n}\sum\limits_{m,s =1}^n Im \hskip 2pt l_{ms}, \ph k= \overline{1, n -1}$;

\noindent
$h_{jk} = 0, \ph j\ne k, \ph j \ne n, \ph k \ne n, \ph j,k = \overline{1,n}$.

\noindent
The matrix $H_L$ is called separator of $L$ and is denoted by $Sep (L)$.

\vskip 5pt
Let $\alpha(t),\ph \beta(t)$ and $\gamma(t)$ be real-valued locally integrable functions on $[t_0,+\infty)$ such that $\alpha(t) + \beta(t) \equiv 1, \ph t \ge t_0$. Set $A_{\alpha,\beta,\gamma}(t) \equiv \alpha(t) A(t) + \beta(t) A^*(t) + \gamma(t) I, \ph t \ge t_0$.
Consider the linear matrix equation
$$
B(t) X = Sep (A_{\alpha,\beta,\gamma}(t)), \phh t \ge t_0. \eqno (3.30)
$$
We are interested in whether this equation has a  Hermitian solution. Indicate some particular cases, when it has a Hermitian solution.

\noindent
IV$^\circ$) $B(t) = \sigma(t) Sep (A_{\alpha,\beta,\gamma}(t)), \ph \sigma(t) \ne 0, \ph t \ge t_0$ Then $H(t) \equiv \frac{1}{\sigma(t)} I$.

\noindent
V$^\circ$) $B(t) = \sigma(t) \sqrt{Sep  (A_{\alpha,\beta,\gamma}(t))}$. Then $H(t) \equiv \frac{1}{\sigma(t)} \sqrt{Sep (A_{\alpha,\beta,\gamma}(t))}$.

\noindent
VI$^\circ$)
$$
Sep  (A_{\alpha,\beta,\gamma}(t))   = diag \biggl\{\stackrel{n_1}{\overbrace{\nu_1(t),\dots,\nu_1(t)}},\dots,\stackrel{n_p}{\overbrace{\nu_p(t),\dots,\nu_p(t)}}, \stackrel{m}{\overbrace{0,\dots,0}}\biggr\}
$$
$n_1 + \dots + n_p + m = n$, $B(t) \equiv Bloc\{B_1(t),\dots,B_p(t),\Theta_m\}$ is a bloc diagonal matrix, $B_k(t) > 0, \ph t \ge t_0$ is a Hermitian matrix of dimension $n_k\times n_k, \ph k=\overline{1,p}, \ph \theta_m$ is a null matrix of dimension $m\times m$. Then $H(t) = Bloc\{\nu_1(t) B_1^{-1}(t),\dots,\nu_p(t) B_p^{-1}(t),\theta_m\}, \ph t \ge t_0$.

\noindent
VII$^\circ$) $Sep(A_{\alpha,\beta,\gamma}(t)) \equiv 0, \ph t \ge t_0.$ Then $H(t)\equiv 0, \ph t \ge t_0$ is a Hermitian solution of Eq.~ (3.30).

 For any $n\times n$ dimensional absolutely continuous matrix function $X(t)$ set
$$
K_X(t) \equiv X'(t) + X(t)B(t)X(t) + A^*(t)X(t) + X(t) A(t) - C(t), \phh t \ge t_0.
$$

{\bf Theorem 3.7.} {\it Let $H(t)$ be an absolutely continuous Hermitian  solution  of Eq. (3.30). If the scalar system
$$
\sist{\phi' =\gamma(t)\phi +  \frac{\lambda_1(B(t))}{n} \phi,}{\psi' = - K_H(t) \phi -\gamma(t)\psi, \ph t \ge t_0}
$$
is oscillatory, then the system (1.1) is also oscillatory.
}

 Proof. In Eq. (2.18) substitute
$$
Y = Z + H(t), \phh t \ge t_0. \eqno (3.31)
$$
We obtain
$$
Z' + Z B(t) Z + (A^*(t) + H(t) B(t)) Z + Z(A(t) + B(t) H(t)) + K_H(t) = 0, \phh t \ge t_1.  \eqno (3.32)
$$
Suppose the system (1.1)  is not oscillatory. Then by (2.19)  and (3.31) Eq. (3.32)  has a solution $Z(t)$ on $[t_1,+\infty)$ for some $t_1\ge t_0$. Hence,
$$
Z'(t) + Z(t) B(t) Z(t) + (A^*(t) + H(t) B(t)) Z(t) + Z(t)(A(t) + B(t) H(t)) + K_H(t) = 0,
$$
$t \ge t_1$. Then since $H(t)$ is a Hermitian solution of Eq. (3.30) we have
$$
Z'(t) + Z(t)B(t)Z(t)+[\alpha(t)(A^*(t)+Sep(A(t))) + \beta(t)(A^*(t)+ Sep(A^*(t)))Y(t)+
$$
$$
Z(t)[\alpha(t)(A(t) + Sep(A(t))) + \beta(t)(A(t) + Sep(A^*(t)))]+K_H(t) + 2\gamma(t) Z(t)= 0,  \eqno (3.33)
$$
$t \ge t_1$.
It was shown in [5], that for any two $n\times n$ dimensional matrices $L$ and $U$ the equality
$$
Sum([L + Sep L] U) = \frac{i \hskip 2pt Im \hskip 2pt (Sum (L))}{n} Sum(U).
$$
is valid. Then
$$
Sum \{(A^*(t) + Sep(A(t))) Z(t) + Z(t)(A(t) + Sep(A(t)))\} = 0.
$$
$$
Sum \{(A^*(t) + Sep(A^*(t))) Z(t) + Z(t)(A(t) + Sep(A^*(t)))\} = 0.
$$
This together with (3.33) and the inequality (see [5])
$$
Sum(Y(t) B(t) Y(t)) \ge \frac{\lambda(B(t))}{n}(Sum(Y(t)))^2, \ph t \ge t_1
$$
implies that
$$
Sum(Z(t))' + \frac{\lambda(B(t))}{n}(Sum(Z(t))^2 + 2\gamma(t)(Sum(Z(t)) +  K_H(t) \le 0, \phh t \ge t_1.
$$
Further as in the proof of Theorem 3.2. The theorem is proved.

{\bf Remark 3.3.} {\it Using (2.23) on the basis of comparison Theorem 2.1 one can show that in Theorem 3.7 $\lambda_1(B(t))$ can be replaced by a more explicit function $\nu_0(B(t))$.
}

\vskip 5pt

{\bf Corollary 2.2.} {\it If

$$
 Sep (A_{\alpha,\beta,\gamma}(t)) \equiv 0, \ph t \ge t_0 \eqno (3.34)
$$
and
$$
\ilp{t_0}\exp\biggl\{-\il{t_0}{t}\gamma(\tau)d\tau\biggr\}\nu_0(B(t)) d t = \ilp{t_0}\exp\biggl\{\il{t_0}{t}\gamma(\tau)d\tau\biggr\}\biggl(- Sum C(t)\biggr) d t = +\infty,
$$
then the system (1.1) is oscillatory.}

Proof. According to IV) if $Sep(A_{\alpha,\beta,\gamma})(t) \equiv 0$ then $H(t)\equiv 0$ is a Hermitian solution of Eq. (3.30)  on $[t_0,+\infty)$. In this case $K_H(t) = - C(t)$. Then Corollary~ 2.2
  follows from  Theorems 2.2, 3.7 and Remark 3.3. The corollary is proved.

\vskip 5pt

{\bf Remark 3.2.} {\it Corollary 2.2 is another generalization of Leighton's oscillation criterion.}

{\bf Example 3.4.} {\it Let $Q_j, \ph j=\overline{1,3}$ be measurable sets such that $Q_j\cap Q_k = \emptyset, \ph j\ne k, \ph j,k = \overline{1,3}$ and $Q_1 \cup Q_2 \cup Q_3 = [t_0,+\infty)$. Let $a_{jk}(t), \ph a_{nk}(t), \ph a_{jn}(t), \ph j,k=\overline{1,n-1}$ be real valued locally integrable functions on $[t_0,+\infty)$. Set
$$
A_{Q_1}(t) \equiv \begin{pmatrix} a_{11}(t) &\dots & a_{1n}(t)\\
\dots & \dots & \dots\\
a_{n-1,1}(t)&\dots & a_{n-1,n}(t)\\
-\sum\limits_{j=1}^{n-1} a_{j1}(t) &\dots & -\sum\limits_{j=1}^{n-1} a_{jn}(t)
\end{pmatrix},
$$
$$
A_{Q_2}(t) \equiv \begin{pmatrix} a_{11}(t) &\dots & a_{1,n-1}(t) & -\sum\limits_{j=1}^{n-1} a_{1j}(t)\\
\dots & \dots & \dots &\dots\\
a_{n,1}(t)&\dots & a_{n,n-1}(t) & -\sum\limits_{j=1}^{n-1} a_{nj}(t)
\end{pmatrix}
$$
$$
A_{Q_3}(t) \equiv \begin{pmatrix} a_{11}(t) &\dots & a_{1,n-1}(t) & -\sum\limits_{j=1}^{n-1} a_{1j}(t)\\
\dots & \dots & \dots &\dots\\
a_{n-1,1}(t)&\dots & a_{n-1,n-1}(t) & -\sum\limits_{j=1}^{n-1} a_{nj}(t)\\
-\sum\limits_{j=1}^{n-1} a_{j1}(t) &\dots & -\sum\limits_{j=1}^{n-1} a_{jn}(t) & \sum\limits_{j,k=1}^{n-1} a_{jk}(t)
\end{pmatrix},
$$
$$
A(t) \equiv \sist{
A_{Q_j}(t), \ph t \in Q_j, \ph j=1,2}{
A_0(t) + A_{Q_3}(t), \ph  t \in Q_{3}}
$$
where $A_0(t)$ is a locally integrable skew symmetric matrix ($A_0^*(t) = - A_0(t)$). Then it is not difficult to verify that for
$
\alpha(t) \equiv \begin{cases}
1, \ph t \in Q_1,\\
0, \ph t \in Q_2,\\
1/2. \ph t \in Q_3
\end{cases}, \ph \gamma(t) \equiv 0, \ph t \ge t_0
$
the condition (3.34) is fulfilled.}

{\bf Example 3.5.} {\it Let $B_0 > 0$ be a Hermitian matrix of dimension $n\times n$. Consider the linear matrix Hamiltonian system
$$
\sist{\Phi' = \frac{t(\sqrt{B_0})^{-1}}{\sqrt{1 + t^2}} \Phi + \frac{B_0}{1+ t^2} \Psi,}{\Psi' = -t \Phi - \frac{t(\sqrt{B_0})^{-1}}{\sqrt{1 + t^2}} \Psi, \ph t \ge t_0.} \eqno (3.35)
$$
It is not difficult to verify that for this system $F(t)\equiv 0$ is a solution of Eq. (3.21) (since for this system $A(t)\sqrt{B(t)} - (\sqrt{B(t)})' \equiv 0, \ph t \ge t_0$). Hence, $J_2(t) = \il{t_0}{t}\frac{tr B_0 d\tau}{1 + T^2} \to +\infty$ for $ t \to +\infty$. By Theorem 3.5 from here it follows that the system (3.25) is oscillatory. Since $\ilp{t_0}\frac{tr B_0}{1+ t^2} d t < + \infty.$ by (2.23) Theorems 3.1 and 3.3 are not applicable to the system (3.35). By (3.18) the condition V) of Theorem 3.4 for the system (3.35) is not fulfilled. Therefore Theorem 3.4 is not applicable to the system (3.35) as well.
}

\vskip 20 pt

\centerline{ \bf References}

\vskip 20pt

\noindent
1. K. I. Al - Dosary, H. Kh. Abdullah and D. Husein. Short note on oscillation of matrix \linebreak \phantom{a} Hamiltonian systems. Yokohama Math. J., vol. 50, 2003.

\noindent
2. R Bellman, Vvedenie v teoriju matric (Russian translation of R. Bellman, Introduction \linebreak \phantom{a}  to matrix algebra. Los Angeles, California, University of southern California), \linebreak \phantom{a} 367 pages.).

\noindent
3. Sh. Chen, Z. Zheng, Oscillation criteria of Yan type for linear Hamiltonian systems, \linebreak \phantom{a} Comput.  Math. with Appli., 46 (2003), 855--862.

\noindent
4.  F. R. Gantmacher, Theory of Matrix. Second Edition (in Russian). Moskow,,\linebreak \phantom{a} ''Nauka'', 1966.

\noindent
5. G. A. Grigorian, Oscillation criteria for linear matrix Hamiltonian systems. \linebreak \phantom{aa} Proc. Amer. Math. Sci, Vol. 148, Num. 8 ,2020, pp. 3407 - 3415.

\noindent
6. G. A. Grigorian.   Interval oscillation criteria for linear matrix Hamiltonian systems,\linebreak \phantom{a}  Rocky Mount. J. Math.,  vol. 50 (2020), No. 6, 2047–2057

\noindent
7.  G. A. Grigorian. Oscillatory criteria for the systems of two first - order Linear \linebreak \phantom{a} ordinary differential equations. Rocky Mount. J. Math., vol. 47, Num. 5,
 2017, \linebreak \phantom{a}  pp. 1497 - 1524

\noindent
8. G. A. Grigorian,  On two comparison tests for second-order linear  ordinary\linebreak \phantom{aa} differential equations (Russian) Differ. Uravn. 47 (2011), no. 9, 1225 - 1240; trans-\linebreak \phantom{aa} lation in Differ. Equ. 47 (2011), no. 9 1237 - 1252, 34C10.

\noindent
9. G. A. Grigorian, Two comparison criteria for scalar Riccati equations and their \linebreak \phantom{aa} applications.
Izv. Vyssh. Uchebn. Zaved. Mat., 2012, Number 11, 20–35.

\noindent
10. G. A. Grigorian, Oscillatory and Non Oscillatory criteria for the systems of two \linebreak \phantom{aa}   linear first order two by two dimensional matrix ordinary differential equations. \linebreak \phantom{aa}   Arch.  Math., Tomus 54 (2018), PP. 189 - 203.

\noindent
11. I. S. Kumary and S. Umamaheswaram, Oscillation criteria for linear matrix \linebreak \phantom{a} Hamiltonian systems, J. Differential Equ., 165, 174--198 (2000).

\noindent
12. A. Kurosh, Higher Algebra, Moskow Mir Publisher (English translation) 1980, \linebreak \phantom{a}428 pages.

\noindent
13. L. Li, F. Meng and Z. Zheng, Oscillation results related to integral averaging technique\linebreak \phantom{a} for linear Hamiltonian systems, Dynamic Systems  Appli. 18 (2009), \ph \linebreak \phantom{a} pp. 725--736.


\noindent
14. Y. G. Sun, New oscillation criteria for linear matrix Hamiltonian systems. J. Math. \linebreak \phantom{a} Anal. Appl., 279 (2003) 651--658.

\noindent
15. C. A. Swanson. Comparison and oscillation theory of linear differential equations.   \linebreak \phantom{a} Academic press. New York and London, 1968.

\noindent
16. Q. Yang, R. Mathsen and S. Zhu, Oscillation theorems for self-adjoint matrix \linebreak \phantom{a}   Hamiltonian
 systems. J. Diff. Equ., 19 (2003), pp. 306--329.

\noindent
17. Z. Zheng, Linear transformation and oscillation criteria for Hamiltonian systems. \linebreak \phantom{a} J. Math. Anal. Appl., 332 (2007) 236--245.

\noindent
18. Z. Zheng and S. Zhu, Hartman type oscillatory criteria for linear matrix Hamiltonian  \linebreak \phantom{a} systems. Dynamic  Systems  Appli., 17 (2008), pp. 85--96.

\end{document}